\documentclass[12pt,a4paper]{article}
\usepackage[cp1251]{inputenc}
\usepackage{amsfonts, amssymb, amsthm, amsmath}
\usepackage{mathtext}
\usepackage[english]{babel}

\usepackage[a4paper, margin=2cm]{geometry}

\begin{document}

\renewcommand{\r}{\mathbb R}
\renewcommand{\le}{\leqslant}
\renewcommand{\ge}{\geqslant}
\newcommand{\eqd}{\stackrel{d}{=}}
\newcommand{\il}[2]{\int\nolimits_{#1}^{#2}}

\title{Max-compound Cox processes. III}

\author{V. Yu. Korolev\textsuperscript{1}, I. A. Sokolov\textsuperscript{2}, A. K.
Gorshenin\textsuperscript{3}}

\date{}

\maketitle

\footnotetext[1]{Faculty of Computational Mathematics and Cybernetics, Lomonosov
Moscow State University; Federal Research Center ``Computer Science and Control'' of the Russian Academy of Sciences; Hangzhou Dianzi University, China; vkorolev@cs.msu.ru}

\footnotetext[2]{Federal Research Center ``Computer Science and Control'' of the Russian Academy of Sciences; Faculty of Computational Mathematics and Cybernetics,
Lomonosov Moscow State University; isokolov@frccsc.ru}

\footnotetext[3]{Federal Research Center ``Computer Science and Control'' of the Russian Academy of Sciences; Faculty of Computational Mathematics and Cybernetics,
Lomonosov Moscow State University; agorshenin@frccsc.ru}

{

\small

{\bf Abstract:} Extreme values are considered in samples with random
size that has a mixed Poisson distribution being generated by a
doubly stochastic Poisson process. We prove some inequalities
providing bounds on the rate of convergence in limit theorems for
the distributions of max-compound Cox processes.

\smallskip


}

\section{Introduction}

In this paper we continue the research we started in
\cite{KorolevSokolovGorshenin2018} and
\cite{KorolevSokolovGorshenin2018b}. The aim of these parers was to
study the analytic properties and asymptotic behavior of
max-compound Cox processes that are defined as extreme order
statistics in samples with random size being doubly stochastic
Poisson (Cox) process. Here we focus on their one-dimensional
distributions, so the presented results can be assigned to the
asymptotic theory of extreme order statistics constructed from
samples with random sizes. The foundation of this theory was laid in
the works \cite{Berman1964, BarndorffNielsen1964, Mogyorodi1967,
GnedenkoSenusiBereksi1982, SenusiBereksiJanic1984}. A review of
these studies can be found in the book \cite{Galambos1984}, also see
\cite{Galambos1994}. In these works the so-called transfer theorems
were proved establishing the convergence of the distributions of
extreme order statistics in samples with random sizes provided that
appropriately normalized extrema in samples with non-random sizes
have limit distributions. For a more recent summary of the results
related to the asymptotic behavior of order statistics with random
sample size see \cite{BakaratShandidy2004a, BakaratShandidy2004b}.

When modeling many real phenomena it can be reasonably assumed that
the flow of informative events each of which generates the next
observation forms a {\it chaotic} point stochastic process on the
time axis. Moreover, this point process can be non-stationary
(time-non-homogeneous) due to stochastic variation of the intensity
caused by external factors. As is known, most reasonable
probabilistic models of non-stationary (time-non-homogeneous)
chaotic point processes are {\it doubly stochastic Poisson
processes} also called {\it Cox processes} (see, e. g.,
\cite{Grandell1976, BeningKorolev2002}). These processes are defined
as Poisson processes with stochastic intensities.

Actually we deal with a particular case of the general limit scheme
considered in the papers mentioned above. However, a special case
considered here, where the sample size is a value of a doubly
stochastic Poisson process, first, generates a rather wide class of
flexible models that can be successfully applied to describe many
real phenomena, and second, allows to easily construct analytic
models with interesting and convenient properties.

In \cite{KorolevSokolovGorshenin2018, KorolevSokolovGorshenin2018b}
we proved limit theorems for the distributions of max-compound Cox
processes. Unlike similar results in \cite{Galambos1984}, in our
theorems the conditions on the tail behavior of the distribution of
observations are separated from the conditions on the asymptotic
behavior of the leading process. In other words, in
\cite{KorolevSokolovGorshenin2018, KorolevSokolovGorshenin2018b} the
conditions for convergence to take place are formulated in terms of
the (cumulative) intensities of the flow of informative events.
Moreover, in the case under consideration it becomes possible to
trace the interrelation of the constants normalizing the extrema
with those normalizing the leading process. The limit distributions
turn out to be power mixtures of the so-called extreme value
distributions of the classical asymptotic extreme value theory.

In \cite{KorolevSokolovGorshenin2018} we proved some transfer
theorems for max-compound Cox processes establishing {\it
sufficient} conditions for the convergence. In
\cite{KorolevSokolovGorshenin2018b} we presented the inverses of
those transfer theorems in order to provide {\it necessary and
sufficient} conditions for the convergence. In the present paper we
focus our attention on proving bounds for the rate of convergence in
limit theorems for max-compound Cox processes.

The paper is organized in the following way. In Section 2 we recall
the definition of a max-compound doubly stochastic Poisson process
and formulate some transfer theorems proved in
\cite{KorolevSokolovGorshenin2018} and
\cite{KorolevSokolovGorshenin2018b}. In Section 3 we apply a new
method of proving convergence rate bounds for extreme order
statistics constructed from samples with the sample size having the
mixed Poisson distribution based on using the classical bounds.
Actually, here we prove some `transfer theorems' for the bounds on
the rate of convergence of max-compound Cox processes. First we
extend the classical bounds to the case where the random sample size
has the binomial distribution. Then using the reasoning based on the
Poisson theorem we transfer these bounds to the case where the
sample size has the Poisson distribution. Finally, the bounds
obtained for the case of the Poisson-distributed sample size are
extended to the case of the random sample size with the mixed
Poisson distribution. This obtained chain of theorems leads to the
conclusion that the case of the random sample size with the Poisson
distribution is in some sense ideal because in this case the
convergence of the distributions of extreme order statistics is more
rapid than in the case with the non-random sample size coinciding
with the parameter of the Poisson distribution. In Section 4 general
bounds of the rate of convergence of the distributions of
max-compound Cox processes are obtained by a direct method.

In what follows the symbols $\Longrightarrow$ and $\eqd$ will denote
convergence in distribution and coincidence of distributions, respectively.
%Euler's gamma function will be denoted $\Gamma(z)$:
%$$
%\Gamma(z)=\il{0}{\infty}x^{z-1}e^{-x}dx,\ \ \ z>0.
%$$

Recall the definition of a doubly stochastic Poisson process (Cox
process).

A {\it doubly stochastic Poisson process} (also called a {\it Cox
process}) is a stochastic point process of the form
$N_1(\Lambda(t))$, where $N_1(t)$, $t\geq0$, is a homogeneous
Poisson process with unit intensity and the stochastic process
$\Lambda(t)$, $t\geq0$, is independent of $N_1(t)$ and possesses the
following properties: $\Lambda(0)=0$, ${\sf P}(\Lambda(t)<\infty)=1$
for any $t>0$, the sample paths of $\Lambda(t)$ do not decrease and
are right-continuous. In this context the Cox process $N(t)$ is said
to be controlled or lead by the process $\Lambda(t)$.

Now let $N(t)$, $t\geq0$, be the a Cox process lead by the process
$\Lambda(t)$. Let $T_1,T_2,\ldots$ be the jump points of the process
$N(t)$. Consider a marked Cox point process $\{(T_i,X_i)\}_{i\ge1}$,
where $X_1,X_2,\ldots$ are independent identically distributed
(i.i.d.) random variables (r.v.'s) independent of the process
$N(t)$. In risk theory, the overwhelming majority of studies related
to the point process $\{(T_i,X_i)\}_{i\ge1}$ deal with {\it compound
Cox process} $S(t)$ that is a function of the marked Cox point
process $\{(T_i,X_i)\}_{i\ge1}$ defined for each $t\ge0$ as the {\it
sum} of all marks of the points of the marked Cox point process
which do not exceed the time $t$. In $S(t)$ the compounding
operation is {\it summation}. Another function of the marked Cox
point process $\{(T_i,X_i)\}_{i\ge1}$ which is often of no less
importance in many applied problems is the so-called max-compound
Cox process that differs from $S(t)$ by that the compounding
summation operation is replaced by the operation of taking the
maximum of the marking r.v.'s.

\smallskip

{\sc Definition 1.} Let $N_{\lambda}(t)$ be a homogeneous Poisson
process with constant intensity $\lambda>0$. The process
$M^{(\lambda)}(t)$ defined as
$$
M^{(\lambda)}(t)=\begin{cases}-\infty, & \text{if $N_{\lambda}(t)=0$,}\vspace{1mm}\\
              {\displaystyle
              \max_{1\le k\le N_{\lambda}(t)}X_k,} & \text{if $N_{\lambda}(t)\ge 1$,}\end{cases}
$$
$t\ge0$, is called a {\it max-compound Poisson process}.

Let $N(t)$ be a Cox process. The process $M(t)$ defined as
$$
M(t)=\begin{cases}-\infty, & \text{if $N(t)=0$,}\vspace{1mm}\\
              {\displaystyle
              \max_{1\le k\le N(t)}X_k,} & \text{if $N(t)\ge 1$,}\end{cases}
$$
$t\ge0$, is called a {\it max-compound Cox process}.

\smallskip

Obviously, a max-compound Poisson process is a special case of a
max-compound Cox process with $\Lambda(t)\equiv\lambda t$.

\smallskip

{\sc Lemma 1} \cite{KorolevSokolovGorshenin2018}. {\it For an
arbitrary $t\ge0$, the d.f. of a max-compound Poisson process has
the form}
$$
{\sf P}\big(M^{(\lambda)}(t)<x\big)=\exp\big\{\lambda
t[F(x)-1]\big\},\ \ \ x\in\r.\eqno(1)
$$

\smallskip

Obviously, Lemma 1 is a particular case of the following more
general proposition.

\smallskip

{\sc Lemma 2} \cite{KorolevSokolovGorshenin2018}. {\it Let
$X_1,X_2,\ldots$ be i.i.d. r.v.'s with the common d.f. $F(x)$. Let
$N$ be an integer-valued nonnegative r.v. with generating function
$\psi(s)$,
$$
\psi(s)={\sf E}s^N=\sum_{n=0}^{\infty}s^n{\sf P}(N=n),\ \ \ 0\le
s\le1.
$$
Assume that the r.v.'s $N,X_1,X_2,\ldots$ are jointly independent.
Then}
$$
{\sf P}\Big(\max_{0\le k\le N}X_k<x\Big)=\psi\big(F(x)\big),\ \ \
x\in\r.
$$

\smallskip

We have $\psi_{\lambda}(s)={\sf E}s^{N_{\lambda}}=e^{\lambda(s-1)}$,
$0\le s\le 1$. Hence Lemma 1 follows from Lemma 2.

\smallskip

Let $M(t)$ be a max-compound Cox process generated by the sequence
$X_1,X_2,\ldots$ and lead by the process $\Lambda(t)$. Then it
follows from Lemma 1 that
$$
{\sf P}\big(M(t)<x\big)={\sf
E}\exp\big\{-\Lambda(t)[1-F(x)]\big\}=\il{0}{\infty}\exp\big\{-\lambda
[1-F(x)]\big\}d{\sf P}\big(\Lambda(t)<\lambda\big),\ \ \ x\in\r.
$$

The relation given above means that in an arbitrary fixed point
$x\in\r$ the d.f. of the max-compound Cox process lead by the
process $\Lambda(t)$ is equal to the Laplace--Stieltjes transform of
the r.v. $\Lambda(t)$ at the point $z=1-F(x)$.

\smallskip

In \cite{KorolevSokolovGorshenin2018} we proved some limit theorems
exposing the asymptotic behavior of max-compound Cox processes and
described the class of limit laws.

The common d.f. of the r.v.'s $X_j$ will be denoted $F(x)$. We will
also use the notation
$$
{\mbox{\rm lext}}(F)=\inf\{x:\ F(x)>0\},\ \ \ {\mbox{\rm
rext}}(F)=\sup\{x:\ F(x)<1\}.
$$

\smallskip

{\sc Theorem 1} \cite{KorolevSokolovGorshenin2018}. {\it Assume that
there exists a positive infinitely increasing function $d(t)$ and a
nonnegative r.v. $\Lambda$ such that
$$
\frac{\Lambda(t)}{d(t)}\Longrightarrow \Lambda\ \ \
(t\to\infty).\eqno(2)
$$
Also assume that $\mbox{\rm rext}(F)=\infty$ and there exists a
positive number $\gamma$ such that
$$
\lim_{y\to\infty}\frac{1-F(yx)}{1-F(y)}=x^{-\gamma}\eqno(3)
$$
for any $x>0$. Then there exists a positive function $b(t)$ and a
d.f. $H_1(x)$ such that
$$
\mbox{\sf P}\Big(\frac{1}{b(t)}\max_{1\le k\le N(t)}X_k<x\Big)
\Longrightarrow H_1(x)\ \ \ (t\to\infty).
$$
Moreover,}
$$
H_1(x)=\il{0}{\infty}e^{-\lambda x^{-\gamma}}d\mbox{\sf
P}(\Lambda<\lambda), \ \ \ x> 0,
$$
{\it and $H_1(x)=0$ for $x\le 0$, whereas the function $b(t)$ can be
defined as}
$$
b(t)=\inf\Big\{x:\ 1-F(x)\le\frac{1}{d(t)}\Big\}.\eqno(4)
$$

\smallskip

{\sc Theorem 2} \cite{KorolevSokolovGorshenin2018}. {\it Assume that
there exist a positive infinitely increasing function $d(t)$ and a
nonnegative r.v. $\Lambda$ such that convergence $(2)$ takes place.
Also assume that $\mbox{\rm rext}(F)<\infty$ and the d.f.
$G_F(x)=F\bigl( \mbox{\rm rext}(F)-x^{-1}\bigr)$ satisfies condition
$(3)$, that is, there exists a positive $\gamma$ such that for any
$x>0$
$$
\lim_{y\to\infty}\frac{1-G_F(yx)}{1-G_F(y)}=x^{-\gamma}.%\eqno(5)
$$
Then there exist functions $a(t)$ and $b(t)$ and a d.f. $H_2(x)$
such that}
$$
\mbox{\sf P}\Big(\frac{1}{b(t)}\Big(\max_{1\le k\le
N(t)}X_k-a(t)\Big)<x\Big) \Longrightarrow H_2(x)\ \ \ (t\to\infty).
$$
{\it Moreover, $H_2(x)=1$ for $x\ge 0$ and}
$$
H_2(x)=\il{0}{\infty}e^{-\lambda|x|^{\gamma}}d\mbox{\sf
P}(\Lambda<\lambda),
$$
{\it for $x<0$. The functions $a(t)$ and $b(t)$ can be defined as}
$$
a(t)={\mbox{\rm rext}}(F),\ \ \ b(t)={\mbox{\rm
rext}}(F)-\inf\Big\{x:\ 1-F(x)\le\frac{1}{d(t)}\Big\}.
$$

\smallskip

{\sc Theorem 3} \cite{KorolevSokolovGorshenin2018}. {\it Assume that
there exist a positive infinitely increasing function $d(t)$ and a
nonnegative r.v. $\Lambda$ such that convergence $(2)$ takes place.
Let
$$
\il{a}{{\mbox{\scriptsize rext}}(F)}[1-F(z)]dz<\infty
$$
for some finite $a$ and for $y\in({\mbox{\rm lext}}(F), {\mbox{\rm
rext}}(F))$ the function $R(y)$ be defined as
$$
R(y)=\frac{1}{1-F(y)}\il{y}{{\mbox{\scriptsize rext}}(F)}[1-F(z)]dz.
$$
Assume that for any $x\in\r$ the limit
$$
\lim_{y\to{\mbox{\scriptsize
rext}}(F)}\frac{1-F\bigl(y+xR(y)\bigr)}{1-F(y)}=e^{-x} %\eqno(5)
$$
exists. Then there exist functions $a(t)$, $b(t)>0$ and $H_3(x)$
such that}
$$
\mbox{\sf P}\Big(\frac{1}{b(t)}\Big(\max_{1\le k\le
N(t)}X_k-a(t)\Big)<x\Big) \Longrightarrow H_3(x)\ \ \ (t\to\infty).
$$
{\it Moreover,}
$$
H_3(x)=\il{0}{\infty}\exp\bigl\{-\lambda e^{-x}\bigr\}d\mbox{\sf
P}(\Lambda< \lambda), \ \ \ x\in\r,
$$
{\it and the functions $a(t)$ can be $b(t)$ chosen in the form}
$$
a(t)=\inf\Big\{x:\ 1-F(x)\le\frac{1}{d(t)}\Big\},\ \ \
b(t)=R(a(t)).%\eqno(5)
$$

\smallskip

The classical asymptotic extreme value theory states that the d.f.'s
$$
H_{1,\gamma}(x)=\begin{cases}\exp\big\{-x^{-\gamma}\big\}, &
\text{if $x> 0$,}\vspace{1mm}\\ 0, & \text{if $x\le 0$,}
\end{cases}
$$
$$
H_{2,\gamma}(x)=\begin{cases}\exp\big\{-|x|^{\gamma}\big\}, &
\text{if $x< 0$,}\vspace{1mm}\\ 1, & \text{if $x\ge 0$,}
\end{cases}
$$
$$
H_{3,0}(x)=\exp\big\{-e^{-x}\big\},\ \ \ x\in\r,
$$
exhaust all possible types of limit laws for normalized maxima of
i.i.d. r.v.'s \cite{Galambos1984, Gumbel1965}. According to von
Mises \cite{Mises1936} and Jenkinson \cite{Jenkinson1955}, these
three types of distributions can be written in a universal form.
Indeed, consider the function
$$
H_{\tau}(x)=\exp\Big\{-\frac{1}{(1+\tau x)^{1/\tau}}\Big\}.
$$
For $\tau=0$ define the function $H_{\tau}(x)$ by its limit value as
$\tau\to0$. Then up to scale and location parameters,
$$
H_{\tau}(x)=\begin{cases}H_{1,\gamma}(x), & \text{for
$\tau>0$,}\vspace{1mm}\\
                   H_{2,\gamma}(x), & \text{for $\tau<0$,}\vspace{1mm}\\
                   H_{3,0}(x), & \text{for $\tau=0$,}\end{cases}
$$
moreover, $\gamma=|1/\tau|$, if $\tau\neq 0$.

\smallskip

{\sc Definition 2.} We shall say that a d.f. $F(x)$ belongs to the
domain of max-attraction of the d.f. $H_{\tau}(x)$, if there exist
sequences of centering and normalizing constants providing the
convergence of the d.f.'s of linearly normalized maxima of i/i/d/
r/v/'s with the common d.f. $F(x)$ to the d.f. $H_{\tau}(x)$.

\smallskip

The fact that the d.f. $F$ belongs to the domain of max-attraction
of the d.f. $H_{\tau}$ will be shortly written as
$F\in\max$-$DA(H_{\tau})$.

The classical asymptotic extreme value theory states that
$F\in\max$-$DA(H_{1, \gamma})$ if and only if $F$ satisfies the
conditions of Theorem 1, $F\in\max$-$DA(H_{2,\gamma})$ if and only
if $F$ satisfies the conditions of Theorem 2,
$F\in\max$-$DA(H_{3,0})$ if and only if $F$ satisfies the conditions
of Theorem 3. Moreover, if a d.f. $F$ does not belong to any of the
domains of max-attraction mentioned above, then the distributions of
extrema in samples from such a distribution cannot converge to any
nondegenerate d.f. whatever centering and normalizing constants are
(e. g., see \cite{Galambos1984}).

\smallskip

{\sc Lemma 3.} {\it If an arbitrary positive function $d(t)$, $t\ge
0$, infinitely increases and $F\in\max$-$DA(H_{\tau})$, then there
exist functions $a(t)$ and $b(t)>0$ such that}
$$
\mbox{\sf P}\!\Big(\frac{1}{b(t)}\Big(\max_{1\le k\le
d(t)}X_k-a(t)\Big)<x\Big) \Longrightarrow H_{\tau}(x)\ \ \
t\to\infty.%\eqno(5)
$$

\smallskip

This statement is proved exactly in the same way as Theorems 1.1 --
1.3 in \cite{BarlowProshan1969}. Moreover, the interrelation of the
functions $a(t)$ and $b(t)$ with $F(x)$ and $d(t)$ is completely
analogous to that described in Theorems 1 -- 3 of
\cite{KorolevSokolovGorshenin2018}.

\smallskip

From Theorems 1 -- 3 we in \cite{KorolevSokolovGorshenin2018} we
obtained the following general form of the transfer theorem for
max-compound Cox processes.

\smallskip

{\sc Theorem 4} \cite{KorolevSokolovGorshenin2018}. {\it Let
$F\in\max$-$DA(H_{\tau})$. Assume that there exist a positive
infinitely increasing function $d(t)$ and a nonnegative r.v.
$\Lambda$ such that convergence $(2)$ takes place. Then there exist
functions $a(t)$ and $b(t)>0$ providing the convergence}
$$
\mbox{\sf P}\!\Big(\frac{1}{b(t)}\Big(\max_{1\le k\le
N(t)}X_k-a(t)\Big)<x\Big) \Longrightarrow
\il{0}{\infty}H_{\tau}^{\lambda}(x)d\mbox{\sf P}(\Lambda<\lambda), \
\ \ t\to\infty.
$$

\smallskip

In the formulations of the theorems presented in Section 2 the
functions $a(t)$, $b(t)$ and $d(t)$ centering and normalizing the
max-compound Cox process and the corresponding leading process were
at least partly specified by the conditions of the theorems. Notice
that the limit mixture distribution in Theorem 4 can be written with
the mixed and mixing distributions determined up to a scale
transform. For example, consider $H_{\tau}(x)$ with $\tau>0$, that
is, let $H_{\tau}(x)=H_{1,\gamma}(x)\equiv
\exp\big\{-x^{-\gamma}\big\}$, $x> 0$, and, correspondingly,
$H_{\tau}(x)=0$ for $x<0$. Then, obviously, the power mixture
$$
H(x)=\int_{0}^{\infty}H_{\tau}^{\lambda}(x)d{\sf
P}(\Lambda<\lambda)=\int_{0}^{\infty}\exp\big\{-\lambda
x^{-\gamma}\big\}d{\sf P}(\Lambda<\lambda),
$$
which is the limit distribution in Theorem 1 and 4, for any
$a\in(0,\infty)$ can also be written as
$$
H(x)=\int_{0}^{\infty}\exp\big\{-\lambda (ax)^{-\gamma}\big\}d{\sf
P}(\Lambda<\lambda d)
$$
with $d=a^{-\gamma}$.

\section{`Transfer theorems' for the estimates of the rate of convergence
in limit theorems for max-compound Cox processes}

In this section we will prove some `transfer theorems' for the
bounds on the rate of convergence of max-compound Cox processes. To
be more precise, we will study how accurate the approximation of the
distributions of max-compound Cox processes by the limit law is.

In \cite{Galambos1984} it was clearly explained that the rate of
convergence of the distributions of extreme order statistics
substantially depends on the centering and normalizing constants and
the following result was proved. Let $H_{\tau}$ be a limit extreme
value distribution and $F\in\max$-$DA(H_{\tau})$. Let $a>0$, $b>0$.
Denote
$$
\zeta_n(x)=n\big(1-F(a+bx)\big),\ \ \ \ \rho_n(x)=\zeta_n(x)+\log
H_{\tau}(x).
$$
It is clear that $\rho_n(x)$ is defined for those $x$, where
$H_{\tau}(x)>0$.

\smallskip

{\sc Theorem 5} \cite{Galambos1984}. {\it Assume that $x$, $a$ and
$b$ satisfy the conditions $H_{\tau}(x)>0$ and $F(a+bx)\ge\frac12$.
Then
$$
\Big|{\sf P}\Big(\max_{1\le i \le
n}X_i<a+bx\Big)-H_{\tau}(x)\Big|\le
H_{\tau}(x)\big[r_{1,n}(x)+r_{2,n}(x)+r_{1,n}(x)r_{2,n}(x)\big],\eqno(5)
$$
where
$$
r_{1,n}(x)=\frac{2\zeta_n^2(x)}{n}+\frac{2\zeta_n^4(x)}{(1-q)n^2}, \
\ \ r_{2,n}(x)=|\rho_n(x)|+\frac{\rho_n^2(x)}{2(1-s)},
$$
and the numbers $q<1$ and $s<1$ satisfy the inequalities
$\frac23\zeta_n^2(x)\le qn$ and $|\rho_n(x)|\le 3s$.}

\smallskip

For the {\sc proof} see \cite{Galambos1984}, Theorem 2.10.1.

\smallskip

Our nearest aim is to transfer the statement of Theorem 5 to the
case of a random sample size with binomial distribution. For this
purpose along with the basic sequence $X_1,X_2,\ldots$ of
independent identically distributed random variables with common
distribution function $F(x)$ consider the numbers $n\in\mathbb{N}$,
$p\in[0,1]$ and a random variable $B_{n,p}$ with the binomial
distribution with parameters $n$ and $p$,
$$
{\sf P}(B_{n,p}=k)=C_n^kp^k(1-p)^{n-k},\ \ \ k=0,1,\ldots,n.
$$
Assume that the random variable $B_{n,p}$ is independent of the
sequence $X_1,X_2,\ldots$. For $i\ge1$ define the random variable
$X_i^*$ as
$$
X_i^*=\begin{cases}\text{\rm lext}(F) & \text{with probability
$1-p$},\cr X_i & \text{with probability $p$}.\end{cases}
$$
For any $x>\text{\rm lext}(F)$ we obviously have
$$
F^*(x)\equiv{\sf P}(X_i^*<x)=1-p+pF(x).%=1-p\big(1-F(x)\big).
$$

\smallskip

{\sc Lemma 4}.
$$
\max_{1\le i\le B_{n,p}}X_i\eqd\max_{1\le i\le n}X_i^*.
$$

\smallskip

{\sc Proof}. By the formula of total probability for any
$x>\text{\rm lext}(F)$ we obviously have
$$
{\sf P}\Big(\max_{1\le i\le
B_{n,p}}X_i<x\Big)=\sum_{k=0}^nC_n^kp^k(1-p)^{n-k}{\sf
P}\Big(\max_{1\le i\le
k}X_i<x\Big)=\sum_{k=0}^nC_n^kp^k(1-p)^{n-k}F^k(x)=
$$
$$
=\sum_{k=0}^nC_n^k\big(pF(x)\big)^k(1-p)^{n-k}=\big(1-p+pF(x)\big)^n=\big(F^*(x)\big)^n=
{\sf P}\Big(\max_{1\le i\le n}X_i^*<x\Big).
$$
The lemma is proved.

\smallskip

Let $a\in\mathbb{R}$, $b>0$. We have
$$
\zeta_n^*(x)\equiv
n\big(1-F^*(a+bx)\big)=n\big(1-1+p-pF(a+bx)\big)=np\big(1-F(a+bx)\big)=p\zeta_n(x).
$$

\smallskip

By virtue of Lemma 4, Theorem 5 implies the following statement.
Denote
$$
\lambda=np.%\eqno(6)
$$

\smallskip

{\sc Theorem 6}. {\it Let the numbers $a$, $b$, $x$, $p$ satisfy the
inequalities $H_{\tau}(x)>0$ and $p\big(1-F(a+bx)\big)\le\frac12$.
Then
$$
\Big|{\sf P}\Big(\max_{1\le i\le
B_{n,p}}X_i<a+bx\Big)-H_{\tau}(x)\Big|\le
H_{\tau}(x)\big[r_{1,n}^*(x)+r_{2,n}^*(x)+r_{1,n}^*(x)r_{2,n}^*(x)\big],
$$
where
$$
r_{1,n}^*(x)=\frac{2\lambda^2\big(1-F(a+bx)\big)^2}{n}+\frac{2\lambda^4\big(1-F(a+bx)\big)^4}{n^2(1-q)},
$$
$$
r_{2,n}^*(x)=\big|\lambda\big(1-F(a+bx)\big)+\log
H_{\tau}(x)\big|+\frac{\big[\lambda\big(1-F(a+bx)\big)+\log
H_{\tau}(x)\big]^2}{2(1-s)},
$$
where $q<1$ and $s<1$ are arbitrary numbers satisfying the
inequalities
$$
\frac{2\lambda^2\big(1-F(a+bx)\big)^2}{3n}\le q%\eqno(6)
$$
and}
$$
\big|\lambda\big(1-F(a+bx)\big)+\log H_{\tau}(x)\big|\le 3s.%\eqno(6)
$$

\smallskip

Now consider the case of a random sample size with the Poisson
distribution. For this purpose along with the basic sequence
$X_1,X_2,\ldots$ of independent identically distributed random
variables with common distribution function $F(x)$ consider a number
$\lambda>0$ and a random variable $N_{\lambda}$ with the Poisson
distribution with parameter $\lambda$,
$$
{\sf P}(N_{\lambda}=k)=e^{-\lambda}\frac{\lambda^k}{k!},\ \ \
k=0,1,\ldots.
$$
Assume that the random variable $N_{\lambda}$ is independent of the
sequence $X_1,X_2,\ldots$. For brevity, we will sometimes use the
notation
$$
\rho_{\lambda}(x)=\lambda\big(1-F(a+bx)\big)+\log H_{\tau}(x).
$$

\smallskip

{\sc Theorem 7}. {\it Let the numbers $a\in\mathbb{R}$, $b>0$,
$x\in\mathbb{R}$, $\lambda>0$, $s\in(0,1)$ satisfy the inequalities
$H_{\tau}(x)>0$ and
$$
\big|\lambda\big(1-F(a+bx)\big)+\log H_{\tau}(x)\big|\le 3s.\eqno(6)
$$
Then}
$$
\Big|{\sf P}\Big(\max_{1\le i\le
N_{\lambda}}X_i<a+bx\Big)-H_{\tau}(x)\Big|\le
H_{\tau}(x)\Big[|\rho_{\lambda}(x)|+\frac{\big(\rho_{\lambda}(x)\big)^2}{2(1-s)}\Big].\eqno(7)
$$

\smallskip

{\sc Proof}. For a given $\lambda>0$ and arbitrary $n\in\mathbb{N}$
consider the random variable $B_{n,\lambda/n}$ having the binomial
distribution with parameters $n$ and $p=\lambda/n$ independent of
the sequence $X_1,X_2,\ldots$. We obviously have
$$
\Big|{\sf P}\Big(\max_{1\le i\le
N_{\lambda}}X_i<a+bx\Big)-H_{\tau}(x)\Big|=\hspace{10cm}
$$
$$
=\Big|{\sf P}\Big(\max_{1\le i\le N_{\lambda}}X_i<a+bx\Big)-{\sf
P}\Big(\max_{1\le i\le B_{n,\lambda/n}}X_i<a+bx\Big)+ {\sf
P}\Big(\max_{1\le i\le
B_{n,\lambda/n}}X_i<a+bx\Big)-H_{\tau}(x)\Big|\le
$$
$$
\le\Big|{\sf P}\Big(\max_{1\le i\le N_{\lambda}}X_i<a+bx\Big)-{\sf
P}\Big(\max_{1\le i\le B_{n,\lambda/n}}X_i<a+bx\Big)\Big|+\Big|{\sf
P}\Big(\max_{1\le i\le
B_{n,\lambda/n}}X_i<a+bx\Big)-H_{\tau}(x)\Big|\equiv
$$
$$
\hspace{12cm}\equiv\Delta_{n,1}(x)+\Delta_{n,2}(x).
$$
First, estimate $\Delta_{n,1}(x)$. Under convention that ${\sf
P}(B_{n,p}=k)=C_n^kp^k(1-p)^{n-k}=0$ for $k>n$, we have
$$
\Delta_{n,1}(x)=\Big|\sum_{k=0}^{\infty}{\sf P}\Big(\max_{1\le i\le
k}X_i<a+bx\Big)\big[{\sf P}(N_{\lambda}=k)-{\sf
P}(B_{n,\lambda/n}=k)\big]\Big|\le\sum_{k=0}^{\infty}\big|{\sf
P}(N_{\lambda}=k)-{\sf P}(B_{n,\lambda/n}=k)\big|.\eqno(8)
$$
To estimate the right-hand side of (8), use the Barbour--Hall
inequality \cite{BarbourHall1984Poisson}, according to which
$$
\sum_{k=0}^{\infty}\big|{\sf P}(N_{\lambda}=k)-{\sf
P}(B_{n,\lambda/n}=k)\big|\le \frac{2\lambda\min\{1,\lambda\}}{n},
$$
and obtain that for a fixed $\lambda=np$ and {\it any}
$n\in\mathbb{N}$ there holds the inequality
$$
\Delta_{n,1}(x)\le\frac{2\lambda\min\{1,\lambda\}}{n}.\eqno(9)
$$
Now estimate $\Delta_{n,2}(x)$ by means of Theorem 6. We have
$$
\Delta_{n,2}(x)\le
H_{\tau}(x)\Big[\frac{2\lambda^2\big(1-F(a+bx)\big)^2}{n}+\frac{2\lambda^4\big(1-F(a+bx)\big)^4}{n^2(1-q)}+
$$
$$
+\big|\lambda\big(1-F(a+bx)\big)+\log H_{\tau}(x)\big|+
\frac{\big[\big|\lambda\big(1-F(a+bx)\big)\big|+\log
H_{\tau}(x)\big]^2}{2(1-s)}+
$$
$$+
\Big(\frac{2\lambda^2\big(1-F(a+bx)\big)^2}{n}+\frac{2\lambda^4\big(1-F(a+bx)\big)^4}{n^2(1-q)}\Big)\times\hspace{8cm}
$$
$$
\hspace{4cm}\times\Big( \big|\lambda\big(1-F(a+bx)\big)+\log
H_{\tau}(x)\big|+\frac{\big[\lambda\big(1-F(a+bx)\big)+\log
H_{\tau}(x)\big]^2}{2(1-s)}\Big)\Big].\eqno(10)
$$
Recall that inequalities (9) and (10) are valid for fixed $\lambda$
and {\it arbitrary} $n$. Let in these inequalities $n\to\infty$. As
this is so, the right-hand side of (9) vanishes whereas (10) turns
into
$$
\Delta_{n,2}(x)\le \big|\lambda\big(1-F(a+bx)\big)+\log
H_{\tau}(x)\big|+
\frac{\big[\big|\lambda\big(1-F(a+bx)\big)\big|+\log
H_{\tau}(x)\big]^2}{2(1-s)},%\eqno(11)
$$
since all the other terms on the right-hand side of (10) vanish as
$n\to\infty$. Moreover, with fixed $a$, $b$, $x$ and $p=\lambda/n$,
as $n\to\infty$, the condition $p\big(1-F(a+bx)\big)\le\frac12$
holds trivially for $n$ large enough. The theorem is proved.

\smallskip

As regards the passage from extrema in samples with
Poisson-distributed sample sizes to those in samples with mixed
Poisson sample sizes, it is connected with some difficulties arizing
because of condition (6). Indeed, the general idea consists in
averaging the right-hand side of (7) to obtain the desired estimate
for mixed Poisson-distributed sample sizes.

Let $N_1(\Lambda)$ be a standard Poisson process considered at a
random point $\Lambda$ that is a nonnegative random variable
independent of the process $N_1(t)$, $t\ge0$. Assume that
$N_1(\Lambda)$ is independent of the sequence $X_1,X_2,\ldots$. As
it follows from Theorems 1 -- 3, the limit law for the extreme order
statistic constructed from the sample $X_1,\ldots,X_{N_1(\Lambda)}$
has the form $H(x)={\sf E}H_{\tau}^{\Lambda}(x)$, so that the
discrepancy to be estimated has the form
$$
\bigg|\int_{0}^{\infty}\Big[{\sf P}\Big(\max_{1\le i\le
N_{\lambda}}X_i<a+bx\Big)-H_{\tau}^{\lambda}(x)\Big]d{\sf
P}(\Lambda<\lambda\Big]\bigg|
$$
and the construction of the reasonable bounds for this discrepancy
by formal averaging of the right-hand side of (7) is difficult due
to that the now the limit law is different being the power mixture
of limit laws for the 'pure' Poisson case. Furthermore, the
necessity to take into account condition (6) makes the everaging
procedure not so trivial.

To avoid all these difficulties we will use another, `direct'
approach to the construction of the bounds for the rate of converge
of max-compound Cox processes (whose distributions are obviosly
mixed Poisson).

\section{General convergence rate estimates in limit theorems for max-compound Cox processes}

Assume that $F\in\max$-$DA(H_\tau)$. For given functions $a(t)$,
$b(t)$ and $d(t)$ and those $x$ which provide the condition
$H_\tau(x)>0$, denote
$$
z_t(x)=d(t)\bigl[1-F\bigl(a(t)+b(t)x\bigr)\bigr],\ \ \
r_t(x)=z_t(x)+\log H_\tau(x).
$$

Show that if the functions $a(t)$, $b(t)$ and $d(t)$ are chosen so
that condition (2) holds and
$$
{\sf P}\!\Big(\frac{1}{b(t)}\Big(\!\max_{1\le k\le N(t)}X_k-a(t)\!
\Big)<x\!\Big)\equiv\il{0}{\infty}\exp\big\{-\lambda z_t(x)\big\}
d{\sf P}(\Lambda(t)<\lambda d(t))\Longrightarrow
%$$
%$$
%\Longrightarrow
\il{0}{\infty}H_\tau^\lambda (x)d{\sf P}(\Lambda<\lambda)\eqno(11)
$$
as $t\to\infty$, then $\sup_tz_t(x)<\infty$ and $r_t(x)\to 0$ as
$t\to\infty$ for $x$ such that $0<H_\tau(x)<1$.

Indeed, by virtue of Theorem 6 in
\cite{KorolevSokolovGorshenin2018b}, from (11) it follows that
$$
{\sf P}\!\Big(\frac{1}{b(t)}\Big(\!\max_{1\le k\le d(t)}X_k-a(t)\!
\Big)<x\!\Big)\Longrightarrow H_\tau(x)\eqno(12)
$$
with the same functions $a(t)$, $b(t)$ and $d(t)$. It is easy to see
that (12) is equivalent to
$$
\exp\big\{-z_t(x)\big\}\to H_\tau(x)\ \ \ (t\to\infty).\eqno(13)
$$
Indeed, let $N_1(t)$ be the standard Poisson process. For any
infinitely increasing function $d(t)$ the relation
$$
\frac{N_1(d(t))}{d(t)}\Longrightarrow 1\ \ \ (t\to\infty)
$$
holds. Therefore, according to Theorem 6 in
\cite{KorolevSokolovGorshenin2018b}, condition (12) is equivalent to
$$
{\sf P}\!\Big(\frac{1}{b(t)}\Big(\!\max_{1\le k\le
N_1(d(t))}X_k-a(t)\! \Big)<x\!\Big)\Longrightarrow H_\tau(x).
$$
But
$$
{\sf P}\!\Big(\frac{1}{b(t)}\Big(\!\max_{1\le k\le
N_1(d(t))}X_k-a(t)\! \Big)<x\!\Big)\equiv\exp\big\{-z_t(x)\big\}.
$$
Thus, the validity of (13) has been established. Passing to
logarithms in (13), we obtain
$$
z_t(x)\to-\log H_\tau(x)\ \ \ (t\to\infty),
$$
but this means that $\sup_tz_t(x)<\infty$ and $r_t(x)\to 0$ as
$t\to\infty$.

The following statement obtained by a direct method without
involvement of the classical estimates for the `non-random' sample
size provides general bound for the convergence rate of the
max-compound Cox process, however, with rather a cumbersome
structure. Below we will consider and discuss some special cases of
this result yielding simpler bounds.

\smallskip

{\sc Theorem 8.} {\it Let functions $a(t)$, $b(t)>0$ and $d(t)>0$ be
arbitrary. For $x$ such that $H_\tau(x)>0$, and any $q\in(0,1)$ and
integer $M\ge1$ there holds the estimate}
$$
\Big|{\sf P}\!\Big(\frac{1}{b(t)}\Big(\!\max_{1\le k\le
N(t)}X_k-a(t)\! \Big)<x\!\Big)-\il{0}{\infty}H_\tau^\lambda(x)d{\sf
P}(\Lambda<\lambda) \Big|\le
$$
$$
\le|r_t(x)|\sum_{k=1}^{M}B_{k,M}(q)|r_t(x)|^{k-1}\il{0}{\infty}\lambda^k
H_\tau^\lambda(x)d{\sf P}(\Lambda<\lambda) +{\sf
P}\!\Big(\Lambda>\frac{q}{|r_t(x)|}\Big)+
$$
$$
+z_t(x)\il{0}{\infty}|{\sf P}(\Lambda(t)<\lambda d(t))- {\sf
P}(\Lambda<\lambda)|e^{-\lambda z_t(x)}d\lambda,%\eqno(14)
$$
{\it where}
$$
B_{k,M}(q)=\begin{cases}{\displaystyle\frac{1}{k!},} &
\text{$k=1,\ldots,M-1$,}\vspace{2mm}\\
{\displaystyle\frac{1+M(1-q)}{(M+1)!(1-q)},} &
\text{$k=M$.}\end{cases}
$$

\smallskip

{\sc Proof}. We have
$$
\Big|{\sf P}\!\Big(\frac{1}{b(t)}\Big(\!\max_{1\le k\le
N(t)}X_k-a(t)\! \Big)<x\!\Big)-\il{0}{\infty}H_\tau^\lambda(x)d{\sf
P}(\Lambda<\lambda) \Big|\le
$$
$$
\le\Big|\il{0}{\infty}e^{-\lambda z_t(x)}d{\sf P}(\Lambda(t)<\lambda
d(t))- \il{0}{\infty}e^{-\lambda z_t(x)}d{\sf
P}(\Lambda<\lambda)\Big|+
$$
$$
+\il{0}{\infty}\Big|e^{-\lambda z_t(x)}-H_\tau^\lambda(x)\Big| d{\sf
P}(\Lambda<\lambda)\equiv I_1(x)+I_2(x).
$$
Consider $I_1(x)$. By integration by parts we make sure that
$$
I_1(x)\le z_t(x)\il{0}{\infty}|{\sf P}(\Lambda(t)<\lambda d(t))-
{\sf P}(\Lambda<\lambda)|e^{-\lambda z_t(x)}d\lambda.\eqno(14)
$$
Consider $I_2(x)$. We have
$$
I_2(x)=\il{0}{q/|r_t(x)|}\Big|e^{-\lambda z_t(x)}-H_\tau^\lambda(x)
\Big|d{\sf P}(\Lambda<\lambda)+
\il{q/|r_t(x)|}{\infty}\Big|e^{-\lambda
z_t(x)}-H_\tau^\lambda(x)\Big| d{\sf P}(\Lambda<\lambda)\equiv
I_{2,1}(x)+I_{2,2}(x).
$$
If $q\in(0,1)$, then
$$
\frac{1}{M!}+\sum_{k=1}^{\infty}\frac{q^k}{(k+M)!}\le
\frac{1}{M!}+\frac{1}{(M+1)!}\sum_{k=1}^{\infty}q^k=\frac{1+M(1-q)}{(M+1)!
(1-q)}.
$$
Therefore rom the Taylor series expansion we can obtain that for
$\lambda|r_t(x)|\le q<1$ there holds the estimate
$$
\Big|e^{-\lambda r_t(x)}-1\Big|\le\sum_{k=1}^{\infty}\frac{\lambda^k
|r_t(x)|^k}{k!}=
%$$
%$$=
|r_t(x)|\Big[\sum_{k=1}^{M-1}\frac{\lambda^k|r_t(x)|^{k-1}}
{k!}+\lambda^M|r_t(x)|^{M-1}\biggl(\frac{1}{M!}+
\sum_{k=1}^{\infty}\frac{q^k}{(k+M)!}\biggr)\Big]\le
$$
$$
\le|r_t(x)|\Big[\sum_{k=1}^{M-1}\frac{\lambda^k|r_t(x)|^{k-1}}
{k!}+\lambda^M|r_t(x)|^{M-1}\frac{1+M(1-q)}{(M+1)!(1-q)}\Big]
$$
Therefore,
$$
I_{21}=\!\!\!\il{0}{q/|r_t(x)|}\!\!\!H_\tau^\lambda(x)
\Big|e^{-\lambda r_t(x)}-1\Big|d{\sf P}(\Lambda<\lambda)\le
$$
$$\le
|r_t(x)|\sum_{k=1}^{M}B_{k,M}(q)\frac{|r_t(x)|^{k-1}}{k!}
\il{0}{q/|r_t(x)|}\lambda^k H_\tau^\lambda(x)d{\sf
P}(\Lambda<\lambda)\le
$$
$$
\le|r_t(x)|\sum_{k=1}^{M}B_{k,M}(q)\frac{|r_t(x)|^{k-1}}{k!}
\il{0}{\infty}\lambda^k H_\tau^\lambda(x)d{\sf P}(\Lambda<\lambda).
\eqno(15)
$$
At the same time
$$
I_{22}\le\!\!\! \il{q/|r_t(x)|}{\infty}\!\!\!d{\sf
P}(\Lambda<\lambda)= {\sf
P}\!\Big(\Lambda>\frac{q}{|r_t(x)|}\Big).\eqno(16)
$$
Now the desired statement follows from (14), (15) and (16). The
theorem is proved.

\smallskip

Setting $M=1$, from Theorem 8 we immediately obtain the following
result.

\smallskip

{\sc Corollary 1.} {\it Under the conditions of Theorem $8$ for any
$q\in(0,1)$ there holds the estimate}
$$
\Big|{\sf P}\!\Big(\frac{1}{b(t)}\Big(\!\max_{1\le k\le
N(t)}X_k-a(t)\! \Big)<x\!\Big)-\il{0}{\infty}H_\tau^\lambda(x)d{\sf
P}(\Lambda<\lambda) \Big|\le
$$
$$
\le\frac{(2-q)|r_t(x)|}{2(1-q)} \il{0}{\infty}\lambda
H_\tau^\lambda(x)d{\sf P}(\Lambda<\lambda)
%\sum_{k=1}^{M}B_{k,M}(q)\frac{|r_t(x)|^{k-1}}{k!}
%\il{0}{\infty}\lambda^k H_\tau^\lambda(x)d{\sf P}(\Lambda<\lambda)
%+\frac{r^2_t(x)(3-2q)}{6(1-q)}\il{0}{\infty}\lambda^2H_\tau^\lambda(x)d{\sf P}
%(\Lambda<\lambda)
+{\sf P}\!\Big(\Lambda>\frac{q}{|r_t(x)|}\Big)+
$$
$$
+z_t(x)\il{0}{\infty}|{\sf P}(\Lambda(t)<\lambda d(t))- {\sf
P}(\Lambda<\lambda)|e^{-\lambda z_t(x)}d\lambda,\eqno(17)
$$

\smallskip

Now assume that $d(t)\equiv t$ and $\Lambda(t)\equiv\lambda t$, that
is, $N(t)\equiv N_{\lambda}(t)$. For this special case Corollary 1
is transformed into the following statement.

\smallskip

{\sc Corollary 2.} {\it assume that $d(t)\equiv t$ and
$\Lambda(t)\equiv\lambda t$, so that $\Lambda =\lambda$, that is,
$N(t)\equiv N_{\lambda}(t)$. There holds the estimate}
$$
\Big|{\sf P}\!\Big(\frac{1}{b(t)}\Big(\!\max_{1\le k\le
N_{\lambda}(t)}X_k-a(t)\! \Big)<x\!\Big)-H_\tau^\lambda(x)
\Big|\le|r_t(x)|\lambda H_\tau^\lambda(x).\eqno(18)
$$

\smallskip

{\sc Example 1}. Let $X_1,X_2,\ldots$ be independent identically
distributed random variables with the common Pareto distribution
function
$$
F(x)=1-\frac{c}{x^{\gamma}+c}, \ \ \ x\ge0,
$$
for some $c>0$, $\gamma>0$. Then
$$
\lim_{y\to\infty}\frac{1-F(xy)}{1-F(y)}=\frac{c(y^{\gamma}+c)}{c(x^{\gamma}y^{\gamma}+c)}=x^{-\gamma},
$$
that is, $F\in DNA(H_{1,\gamma})$,
$H_{1,\gamma}(x)=e^{-x^{-\gamma}}$. Assume that $\Lambda(t)\equiv
\Lambda t$, $t\ge0$, where $\Lambda$ is a random variable with the
gamma distribution given by its probability density
$$
g(x;r,1)=\frac{1}{\Gamma(r)}x^{r-1}e^{-x},\ \ \ x\ge0,
$$
with some $r>0$ (see Example 5 in
\cite{KorolevSokolovGorshenin2018}). Then the distribution of $N(t)$
is negative binomial with parameters $r$ and $p=(t+1)^{-1}$: for
$k=0,1,2,\ldots$
$$
{\sf
P}\big(N(t)=k\big)=\frac{t^k}{k!\Gamma(r)}\int_{0}^{\infty}\lambda^{k+r-1}e^{-\lambda(t+1)}d\lambda=
\frac{t^k}{k!\Gamma(r)(t+1)^{k+r}}\int_{0}^{\infty}\lambda^{k+r-1}e^{-\lambda}d\lambda=
$$
$$
=\frac{\Gamma(k+r)}{k!\Gamma(r)}\Big(\frac{1}{t+1}\Big)^r\Big(1-
\frac{1}{t+1}\Big)^k.
$$
Correspondingly, $t=(1-p)p^{-1}$. Obviously,
$t^{-1}\Lambda(t)\Longrightarrow \Lambda$ as $t\to\infty$, so we can
set $d(t)\equiv t$ and $a(t)\equiv 0$. Let $N_{r,p}$ denote a random
variable having the negative binomial distribution with parameters
$r>0$ and $p\in(0,1)$ independent of the sequence $X_1,X_2,\ldots$.
Now consider the rate of convergence, as $p\to0$, of the
distribution of the normalized maxima
$$
M_{r,p}=\frac{p}{1-p}\max_{1\le k\le N_{r,p}}X_k
$$
to the limit law
$$
H_1(x)=\frac{1}{\Gamma(r)}\int_{0}^{\infty}H_{1,\gamma}^{\lambda}\lambda^{r-1}e^{-\lambda}d\lambda=
\frac{1}{\Gamma(r)}\int_{0}^{\infty}e^{-\lambda
x^{-\gamma}}\lambda^{r-1}e^{-\lambda}d\lambda=\Big(\frac{x^{\gamma}}{1+x^{\gamma}}\Big)^r,\
\ \ \ x\ge0,
$$
prescribed by Theorem 1. In accordance with (18) we have
$$
\Big|{\sf
P}(M_{r,p}<x)-\Big(\frac{x^{\gamma}}{1+x^{\gamma}}\Big)^r\Big|=
$$
$$
=\bigg|\frac{1}{\Gamma(r)}\int_{0}^{\infty}{\sf P}\Big(\frac1t
\max_{1\le k\le
N_{\lambda}(t)}X_k<x\Big)\lambda^{r-1}e^{-\lambda}d\lambda-\frac{1}{\Gamma(r)}\int_{0}^{\infty}e^{-\lambda(x^{-\gamma}+1)}\lambda^{r-1}d\lambda\bigg|\le
$$
$$
\le\frac{1}{\Gamma(r)}\int_{0}^{\infty}\Big|{\sf P}\Big(\frac1t
\max_{1\le k\le
N_{\lambda}(t)}X_k<x\Big)-H_{1,\gamma}^{\lambda}(x)\Big|\lambda^{r-1}e^{-\lambda}d\lambda\le
\frac{1}{\Gamma(r)}\int_{0}^{\infty}|r_t(x)|
H_{1,\gamma}^{\lambda}(x)\lambda^re^{-\lambda}d\lambda=
$$
$$
=r|r_t(x)|\Big(\frac{x^{\gamma}}{1+x^{\gamma}}\Big)^{r+1}.
$$
Moreover, in this case for $t>1$ we can take
$b(t)=F^{-1}(1-\frac{1}{t})$, so we have
$$
|r_t(x)|=t\Big[1-F\Big(xF^{-1}\Big(1-\frac{1}{t}\Big)\Big)\Big]+\log
H_{1,\gamma}(x)=\frac{t}{x^{\gamma}(t-1)+1}-\frac{1}{x^{\gamma}}\le\frac{1}{x^{\gamma}(t-1)},
$$
and recalling that $t=(1-p)p^{-1}$, we finally obtain for
$0<p\le\frac12$:
$$
\sup_{x\ge0}\Big|{\sf
P}(M_{r,p}<x)-\Big(\frac{x^{\gamma}}{1+x^{\gamma}}\Big)^r\Big|\le
\frac{pr}{1-2p}\cdot\sup_{y\ge0}\frac{y^r}{(1+y)^{r+1}}=\frac{p}{1-2p}\cdot\Big(\frac{r}{r+1}\Big)^r.
$$
That is, in this case, as $p\to0$, the uniform distance is $O(p)$.

\smallskip

Since $N_{\lambda}(t)\eqd N_t(\lambda)\eqd N_{\lambda t}(1)\eqd
N_1(\lambda t)$, in (18) the parameters $t$ and $\lambda$ are in
some sense interchangeable, so that in (18) we can set $t=\lambda$,
$a(\lambda)=a$, $b(\lambda)=b$ $d(\lambda)\equiv\lambda$ so that
$\Lambda=1$ and the following inequality for the distribution of
extreme order statistic in a sample with random size having the
Poisson distribution holds.

\smallskip

{\sc Corollary 3.} {\it There holds the estimate}
$$
\big|{\sf P}\big(\max_{1\le k\le
N_{\lambda}}X_k<a+bx\big)-H_{\tau}(x)\big|\le|\rho_{\lambda}(x)|
H_\tau(x),\eqno(19)
$$
{\it where $\rho_{\lambda}(x)=\lambda\big(1-F(a+bx)\big)+\log
H_{\tau}(x)$.}

\smallskip

We see that, first, the right-hand side of (18) is much simpler than
that of (7) and, second, here we do not use any conditions relating
$a$, $b$, $\lambda$ and $x$ like (6).

\smallskip

{\sc Remark 1}. Comparing the structure of the bounds for
convergence rate obtained in Theorems 5 with those established by
Theorem 7 and Corollary 3, we see that the right-hand sides of (7)
and (19) are considerably simpler than that of (5). This may mean
that the case of a Poisson-distributed sample size or, which is the
same, of max-compound Poisson processes is a kind of an ideal model
for which the distributions of extreme order statistics converge
{\it to the same limit law} much faster, than those in the
`classical' case of non-random sample size. This leads to the
conclusion that, like ordinary compound Poisson distributions play
an important role in the classical summation theory providing the
method of accompanying infinitely divisible distributions,
max-compound Poisson processes can play a significant role in the
classical extreme value theory.

\smallskip

Return to the general case. Using the Markov inequality to estimate
the second term on the right-hand side of (17) and minimizing the
result with respect to $q\in(0,1)$, we obtain the following
statement.

\smallskip

{\sc Corollary 4.} {\it Assume that ${\sf E}\Lambda=L<\infty$. Then}
$$
\Big|{\sf P}\!\Big(\frac{1}{b(t)}\Big(\!\max_{1\le k\le
N(t)}X_k-a(t)\! \Big)<x\!\Big)-\il{0}{\infty}H_\tau^\lambda(x)d{\sf
P}(\Lambda<\lambda) \Big|\le
$$
$$
\le|r_t(x)|\Big[\il{0}{\infty}\lambda H_\tau^\lambda(x)d{\sf P}
(\Lambda<\lambda)+L+\Big(2L\il{0}{\infty}\lambda
H_\tau^\lambda(x)d{\sf P} (\Lambda<\lambda)\Big)^{1/2}\Big]+
$$
$$
+z_t(x)\il{0}{\infty}|{\sf P}(\Lambda(t)<\lambda d(t))- {\sf
P}(\Lambda<\lambda)|e^{-\lambda z_t(x)}d\lambda.
$$

\smallskip

Since $0\le H_\tau(x)\le1$, then $H_\tau^\lambda(x)\le1$ for any
$\lambda\ge0$. Hence, for any $k\ge1$
$$
\il{0}{\infty}\lambda^k H_\tau^\lambda(x)d{\sf
P}(\Lambda<\lambda)\le \il{0}{\infty}\lambda^kd{\sf
P}(\Lambda<\lambda)={\sf E}\Lambda^k.
$$
Therefore, Corollary 4 immediately implies the following resiult.

\smallskip

{\sc Corollary 5}. {\it Assume that ${\sf E}\Lambda=L<\infty$. Then}
$$
\Big|{\sf P}\!\Big(\frac{1}{b(t)}\Big(\!\max_{1\le k\le
N(t)}X_k-a(t)\! \Big)<x\!\Big)-\il{0}{\infty}H_\tau^\lambda(x)d{\sf
P}(\Lambda<\lambda) \Big|\le
$$
$$
\le(2+\sqrt{2})L|r_t(x)|
%$$
%$$
+z_t(x)\il{0}{\infty}|{\sf P}(\Lambda(t)<\lambda d(t))- {\sf
P}(\Lambda<\lambda)|e^{-\lambda z_t(x)}d\lambda.
$$

\smallskip

It is easy to make sure that for $\alpha\in(0,1)$ and $b>0$ the
relation
$$
\sup_{\lambda\ge0}\lambda^b\alpha^{\lambda}=\biggl(\frac{b}{e\log
(1/\alpha)}\biggr)^b
$$
holds, Therefore, for any $k\ge1$
$$
\il{0}{\infty}\lambda^k H_\tau^\lambda(x)d{\sf
P}(\Lambda<\lambda)\le \biggl(\frac{k}{e\log(1/\alpha)}\biggr)^k.
$$
Hence, Corollaries 1 and 2, in turn, imply the following statement
that is meaningful even for for ${\sf E}\Lambda=\infty$.

\smallskip

{\sc Corollary 6.} {\it Under the conditions of Theorem $8$, for any
$q\in(0,1)$ there holds the estimate}
$$
\Big|{\sf P}\!\Big(\frac{1}{b(t)}\Big(\!\max_{1\le k\le
N(t)}X_k-a(t)\! \Big)<x\!\Big)-\il{0}{\infty}H_\tau^\lambda(x)d{\sf
P}(\Lambda<\lambda) \Big|\le
$$
$$
\le\frac{(2-q)|r_t(x)|}{2e(1-q)\log\bigl(1/H_\tau(x)\bigr)}
%\il{0}{\infty}\lambda H_\tau^\lambda(x)d{\sf P}(\Lambda<\lambda)
%\sum_{k=1}^{M}B_{k,M}(q)\frac{|r_t(x)|^{k-1}}{k!}
%\il{0}{\infty}\lambda^k H_\tau^\lambda(x)d{\sf P}(\Lambda<\lambda)
%+\frac{r^2_t(x)(3-2q)}{6(1-q)}\il{0}{\infty}\lambda^2H_\tau^\lambda(x)d{\sf P}
%(\Lambda<\lambda)
+{\sf P}\!\Big(\Lambda>\frac{q}{|r_t(x)|}\Big)+
%$$
%$$+
z_t(x)\il{0}{\infty}|{\sf P}(\Lambda(t)<\lambda d(t))- {\sf
P}(\Lambda<\lambda)|e^{-\lambda z_t(x)}d\lambda,
$$

\smallskip

{\sc Example 2.} Assume that $F(x)=(1-e^{-x}){\bf 1}(x\ge0)$ and
$\Lambda (t)\equiv t$ (this is equivalent to that $N(t)\equiv
N_1(t)$). Then $\Lambda(t)/t\equiv 1$, so that the function $d(t)$
can be reasonably determined as $d(t)\equiv t$. In this case
$$
{\sf P}\!\Big(\frac{1}{b(t)}\Big(\!\max_{1\le k\le N(t)}X_k-a(t)\!
\Big)<x\!\Big)=\exp\big\{-te^{-b(t)x-a(t)}\big\}.
$$
Now if we take $b(t)\equiv 1$ and $a(t)\equiv\log t$, then we obtain
$$
{\sf P}\!\Big(\frac{1}{b(t)}\Big(\!\max_{1\le k\le N(t)}X_k-a(t)\!
\Big)<x\!\Big)\equiv\exp\big\{-e^{-x}\big\},\ \ \ x\in\r,
$$
that is, under this choice of the functions $a(t)$, $b(t)$ and
$d(t)$ we have
$$
{\sf P}\!\Big(\frac{1}{b(t)}\Big(\!\max_{1\le k\le N(t)}X_k-a(t)\!
\Big)<x\!\Big)\equiv H_{3,0}(x),\ \ \ x\in\r,
$$
without any passage to the limit. In the case under consideration
$z_t(x)=e^{-x}$ and $r_t(x)\equiv 0$. Therefore, from Corollary 4 we
obtain the bound
$$
\Big|{\sf P}\!\Big(\max_{1\le k\le N(t)}X_k-\log t\!
<x\!\Big)-H_{3,0}(x)\Big|\le 0.
$$

\smallskip

\subsection*{Acknowledgement} The results of Sections 2 and 3 were
formulated and proved by V. Korolev and A. Gorshenin who were
supported by the Russian Science Foundation, project 18-11-00155.

\end{document}